\documentclass[a4paper,12pt]{amsart}
\usepackage{amsfonts}
\usepackage{amssymb}
\usepackage{ifthen}
\usepackage{graphicx}
\nonstopmode \numberwithin{equation}{section}
\usepackage[usenames]{color}
\newtheorem{thm}{Theorem}[section]
\newtheorem{lem}{Lemma}[section]
\newtheorem{cor}{Corollary}[section]
\newtheorem{prop}{Proposition}[section]

\newtheorem{cl}{Claim}[section]
\newtheorem{ca}{Case}[section]
\newtheorem{sca}{Subcase}[section]
\newtheorem{scl}{Subclaim}[section]
\newtheorem{conj}[equation]{Conjecture}

\theoremstyle{definition}
\newtheorem{defn}{Definition}[section]
\newtheorem{op}[equation]{Open Problem}
\newtheorem{ques}[equation]{Question}
\newtheorem{rem}{Remark}[section]
\newtheorem{exam}[equation]{Example}

\newcounter {own}
\def\theown {\thesection       .\arabic{own}}

\newenvironment{pf}[1][]{%
 \vskip 3mm
 \noindent
 \ifthenelse{\equal{#1}{}}%
  {{\slshape Proof. }}%
  {{\slshape #1.} }%
 }%
{\qed\bigskip}

\newcounter{alphabet}
\newcounter{tmp}
\newenvironment{Thm}[1][]{\refstepcounter{alphabet}%
\bigskip%
\noindent%
{\bf Theorem \Alph{alphabet}}%
\ifthenelse{\equal{#1}{}}{}{ (#1)}%
{\bf .} \itshape}{\vskip 8pt}

\makeatletter
\newcommand{\Ref}[1]{\@ifundefined{r@#1}{}{\setcounter{tmp}{\ref{#1}}\Alph{tmp}}}
\makeatother

\newcommand{\IC}{{\mathbb C}}
\newcommand{\ID}{{\mathbb D}}

\def\be{\begin{equation}}
\def\ee{\end{equation}}

\newcommand{\ben}{\begin{enumerate}}
\newcommand{\een}{\end{enumerate}}

\newcommand{\blem}{\begin{lem}}
\newcommand{\elem}{\end{lem}}
\newcommand{\bthm}{\begin{thm}}
\newcommand{\ethm}{\end{thm}}
\newcommand{\bcor}{\begin{cor}}
\newcommand{\ecor}{\end{cor}}
\newcommand{\beg}{\begin{exam}}
\newcommand{\eeg}{\end{exam}}
\newcommand{\begs}{\begin{examples}}
\newcommand{\eegs}{\end{examples}}
\newcommand{\bdefe}{\begin{defn}}
\newcommand{\edefe}{\end{defn}}
\newcommand{\bprob}{\begin{prob}}
\newcommand{\eprob}{\end{prob}}
\newcommand{\bques}{\begin{ques}}
\newcommand{\eques}{\end{ques}}
\newcommand{\bei}{\begin{itemize}}
\newcommand{\eei}{\end{itemize}}
\newcommand{\bcon}{\begin{conj}}
\newcommand{\econ}{\end{conj}}
\newcommand{\bop}{\begin{op}}
\newcommand{\eop}{\end{op}}

\newcommand{\bca}{\begin{ca}}
\newcommand{\eca}{\end{ca}}
\newcommand{\bsca}{\begin{sca}}
\newcommand{\esca}{\end{sca}}

\newcommand{\bcl}{\begin{cl}}
\newcommand{\ecl}{\end{cl}}

\newcommand{\bscl}{\begin{scl}}
\newcommand{\escl}{\end{scl}}

\newcommand{\bcons}{\begin{conjs}}
\newcommand{\econs}{\end{conjs}}

\newcommand{\bprop}{\begin{prop}}
\newcommand{\eprop}{\end{prop}}

\newcommand{\br}{\begin{rem}}
\newcommand{\er}{\end{rem}}
\newcommand{\brs}{\begin{rems}}
\newcommand{\ers}{\end{rems}}
\newcommand{\bo}{\begin{obser}}
\newcommand{\eo}{\end{obser}}
\newcommand{\bos}{\begin{obsers}}
\newcommand{\eos}{\end{obsers}}
\newcommand{\bpf}{\begin{pf}}
\newcommand{\epf}{\end{pf}}
\newcommand{\ba}{\begin{array}}
\newcommand{\ea}{\end{array}}
\newcommand{\beq}{\begin{eqnarray}}
\newcommand{\beqq}{\begin{eqnarray*}}
\newcommand{\eeq}{\end{eqnarray}}
\newcommand{\eeqq}{\end{eqnarray*}}

\newcounter{minutes}
\divide\time by 60
\newcounter{hours}
\multiply\time by 60 \addtocounter{minutes}{-\time}

\begin{document}

\bibliographystyle{amsplain}
\title{Several properties of $\alpha$-harmonic functions in the unit disk}

\def\thefootnote{}
\footnotetext{ \texttt{\tiny File:~\jobname .tex,
          printed: \number\year-\number\month-\number\day,
          \thehours.\ifnum\theminutes<10{0}\fi\theminutes}
} \makeatletter\def\thefootnote{\@arabic\c@footnote}\makeatother

\author{Peijin Li}
\address{Peijin Li,
Department of Mathematics,
Shantou University, Shantou, Guangdong 515063, People's Republic of China and Department of Mathematics,
Hunan First Normal University, Changsha, Hunan 410205, People's Republic of China}
\email{wokeyi99@163.com}

\author{Xiantao Wang$^{~\mathbf{*}}$}
\address{Xiantao Wang, Department of Mathematics,
Shantou University, Shantou, Guangdong 515063, People's Republic of China}
\email{xtwang@stu.edu.cn}

\author{Qianhong Xiao}
\address{Qianhong Xiao, Department of Mathematics,
Shantou University, Shantou, Guangdong 515063, People's Republic of China}
\email{14qhxiao@stu.edu.cn}

\date{}
\keywords{Laplace differential operator, Coefficient estimate, Schwarz-Pick estimate, Landau type theorem.\\
$^{\mathbf{*}}$Corresponding author}

\begin{abstract}
The aim of this paper is to obtain the Schwarz-Pick type inequality for $\alpha$-harmonic functions $f$ in the unit disk and get estimates on the coefficients of $f$. As an application, a Landau type theorem of $\alpha$-harmonic functions is established.
\end{abstract}

\maketitle \pagestyle{myheadings} \markboth{ Peijin Li, Xiantao Wang and Qianhong Xiao}{Several properties of $\alpha$-harmonic functions in the unit disk}

\section{Introduction and main results}\label{sec-1}
Let $\mathbb{C}$ be the complex plane. For $a\in$ $\mathbb{C}$, let $\mathbb {D}(a,r)=\{z:|z-a|<r\}$ $(r>0)$  and $\mathbb {D}(0,r)=\mathbb {D}_{r}$. Also, we use the notations $\mathbb {D}=\mathbb {D}_1$ and $\mathbb{T}=\partial\mathbb{D}$, the boundary of $\mathbb{D}$.

Let
$$
A=\left(
  \begin{array}{cc}
     a & b \\
     c & d \\
  \end{array}
\right)\in \mathbb{R}^{2\times2}.
$$
We will consider the matrix norm
$$|A|=\sup\{|Az|:\; z\in \mathbb{C},\;|z|=1\}$$
and the matrix function
$$l(A)=\inf\{|Az|:\;z\in \mathbb{C},\;|z|=1\}.$$

Let $D$ and $\Omega $ be domains in $\IC$, and let $f=u+iv$: $D\to \Omega$ be a function
that has both partial derivatives at $z=x+iy$ in $D$, where $u$ and $v$ are real functions. The Jacobian matrix of $f$ at $z$ is denoted by
$$
D f(z) =
\left(
  \begin{array}{cc}
    u_x & u_y \\
    v_x & v_y \\
  \end{array}
\right).$$

Set $$\frac{\partial }{\partial{z}}=\frac{1}{2}\Big(\frac{\partial }{\partial{x}}-i\frac{\partial }{\partial{y}}\Big)\;\;\mbox{and}\;\;
\frac{\partial }{\partial{\overline{z}}}=\frac{1}{2}\Big(\frac{\partial }{\partial{x}}+i\frac{\partial }{\partial{y}}\Big).$$
Then
\be\label{qh-1}
|D f(z)|=\sup\{|D f(z)\varsigma|:\; |\varsigma|=1\}=|f_z(z)|+|f_{\overline z}(z)|,
\ee
\be\label{qh-2}
l(D f(z))=\inf\{|D f(z)\varsigma|:\; |\varsigma|=1\}=\big ||f_z(z)|-|f_{\overline z}(z)|\big|
\ee
and
$$|J_{f}(z)|=|D f(z)|\cdot l(D f(z)),$$ where $J_{f}(z)$ stands for the Jacobian of $f$ at $z$.

We denote by $\Delta _{\alpha}$ the weighted Laplace operator corresponding to the so-called standard weight $w_{\alpha}=(1-|z|^{2})^{\alpha}$, that is,
$$
\Delta_{\alpha ,z}=\frac{\partial}{\partial z}(w_{\alpha})^{-1}\frac{\partial}{\partial \overline{z}}=\frac{\partial}{\partial z}(1-{|{z}|}^{2})^{-\alpha}\frac{\partial}{\partial \overline{z}}
$$
in $\mathbb{D}$, where $\alpha > -1$ (see \cite[Proposition $1.5$]{oA} for the reason for this constraint).

In \cite{oA}, Olofsson and Wittsten introduced this operator $\Delta _{\alpha}$ and a counterpart of the classical Poisson integral formula was given.

We remark that in the study of Bergman spaces of $\mathbb{D}$, one often considers the weights $w_{\alpha}$
in $\mathbb{D}$ $(\alpha>-1)$. For an account of recent developments in
Bergman space theory, we mention the monograph \cite{HKZ} by Hedenmalm, Korenblum and Zhu. The case $\alpha= 0$ is commonly referred to as the unweighted case,
whereas the case $\alpha= 1$ has attracted special attention recently with contributions
by Hedenmalm, Shimorin and others (see for instance \cite{HeO, HeP, HeS, SSh} etc).

 Of particular interest to us is the following $\alpha$-harmonic equation in $\mathbb{D}$:
\beq\label{eq1.1}
\Delta_{\alpha}(f)=0.
\eeq

Denote the associated \emph{Dirichlet boundary value problem} of functions $f$ satisfying the equation \eqref{eq1.1} by
\beq\label{eq1.2}
\left\{
\begin{aligned}
\Delta_{\alpha}(f)&=0\;\;\;\;\text{in}\;\mathbb{D},\\
                 f&=f^{\ast}\;\;\text{on}\;\mathbb{T}.
\end{aligned}
\right.
\eeq
Here the boundary data $f^{\ast}$ is a distribution on $\mathbb{T}$, i.e. $f^{\ast}\in \mathcal{D}^{\prime}(\mathbb{T})$, and the boundary condition in the equation \eqref{eq1.2} is to be understood as $f_{r}\to f^{\ast}\in \mathcal{D}^{\prime}(\mathbb{T})$ as $r \to 1^{-}$, where

\be\label{eq1.3}
f_{r}(e^{i\theta})=f(re^{i\theta})
\ee
for $\theta\in [0, 2\pi]$ and $r\in[0,1)$.

For simplicity, we introduce the following definition.
\bdefe
For $\alpha>-1$, a complex-valued function $f$ is said to be {\it $\alpha$-harmonic} if $f$ is twice continuously differentiable in $\mathbb{D}$  and satisfies the condition \eqref{eq1.1}.
\edefe

In \cite{oA}, Olofsson and Wittsten showed that if an $\alpha$-harmonic function $f$ satisfies
$$\lim_{r\to 1^{-}}f_{r}=f^{\ast} \in \mathcal{D}^{\prime}(\mathbb{T})\;\; (\alpha > -1),$$ then it has the form of a \emph{Poisson type integral}
\beq\label{eq1.4}\;\;\;
f(z)=\frac{1}{2\pi}\int^{2\pi}_{0}\mathcal{P}_{\alpha}(ze^{-i\theta})f^{\ast}(e^{i\theta})d\theta
\eeq in $\mathbb{D}$,
where
$$\mathcal{P}_{\alpha}(z)=\frac{(1-|z|^{2})^{\alpha+1}}{(1-z)(1-\overline{z})^{\alpha+1}}.$$

In the following, we always assume that any $\alpha$-harmonic function has such a representation which plays a key role in the discussions of this paper.

Obviously, $\alpha$-harmonicity coincides with harmonicity when $\alpha = 0$. See \cite{DuR} and the references therein for the properties of harmonic mappings.
Particularly, Colonna proved the following Schwarz-Pick type inequality.

\begin{Thm}\label{Thm C}  $($\cite[Theorems 3 and 4]{F}$)$
Let $f$ be a harmonic function of $\mathbb{D}$ into $\mathbb{D}$. Then for $z\in\mathbb{D}$,
$$
|D{f}(z)|\leq\frac{4}{\pi}\cdot\frac{1}{1-|z|^{2}}.
$$
This estimate is sharp and all the extremal functions are
$$
f(z)=\frac{2\delta}{\pi}\arg\Big(\frac{1+\psi(z)}{1-\psi(z)}\Big),
$$
where $\delta\in \mathbb{C}$, $|\delta|=1$ and $\psi$ is a conformal automorphism of $\mathbb{D}$.
\end{Thm}

For the related discussions on this topic, see \cite{HH, cR, cSW, CM, KaM, M} etc.

As the first aim of this paper, we shall generalize Theorem \Ref{Thm C} to the case of $\alpha$-harmonic functions. Our first result is as follows.

\begin{thm}\label{thm1.1}
Suppose that $f$ is an $\alpha$-harmonic function in $\mathbb{D}$ with $\alpha>-1$, that $f^{*}\in C(\mathbb{T})$ and that $\sup_{z\in\overline{\mathbb{D}}}|f(z)|\leq M$, where $M$ is a constant. Then
for $z\in\mathbb{D}$,
$$
|D{f}(z)|\leq\frac{M(\alpha+2)}{\mathrm{c}_{\alpha}}\cdot\frac{1}{1-|z|}\leq\frac{2M(\alpha+2)}{\mathrm{c}_{\alpha}}\cdot\frac{1}{1-|z|^2},
$$
where $\mathrm{c}_{\alpha}=\frac{(\Gamma(\frac{\alpha}{2}+1))^{2}}{\Gamma(\alpha+1)}$ and  $\Gamma(s)=\int^{\infty}_{0}t^{s-1}e^{-t}dt$ $(s>0)$ is the Gamma function.

In particular, if $f$ maps $\ID$ into $\ID$, then
$$|D{f}(z)|\leq\frac{2(\alpha+2)}{\mathrm{c}_{\alpha}}\cdot\frac{1}{1-|z|^2}.$$
\end{thm}

Let $\lambda_{D}(z)|dz|$ be the hyperbolic metric of the domain $D$ having constant Gaussian curvature $-1$. The hyperbolic distance $d_{h_D}(z_1, z_2)$ between two points $z_1$ and $z_2$ in $D$ is defined by
$$
\inf_{\gamma}\Big\{\int_{\gamma}\lambda_{D}(z)|dz|\Big\},
$$
where the infimum is taken over all rectifiable curves $\gamma$ in $D$ connecting $z_1$ and $z_2$.

We have known that if $D=\mathbb{D}$, then (cf. \cite{Bea})
$$\lambda_{\ID}(z)=\frac{2}{1-|z|^2}\;\;\mbox{and}\;\;d_{h_\mathbb{D}}(z_1, z_2)=\log\frac{|1-z_1\overline{z}_2|+|z_1-z_2|}{|1-z_1\overline{z}_2|-|z_1-z_2|}.$$

As a consequence of Theorem \ref{thm1.1}, we have
\bcor\label{cor1.1}
Under the assumptions of Theorem \ref{thm1.1}, if $f$ maps $\ID$ into $\ID$, then
for $z_1$ and $z_2\in \ID$,
$$|f(z_1)-f(z_2)|\leq \frac{\alpha+2}{\mathrm{c}_{\alpha}}d_{h_\mathbb{D}}(z_1, z_2).$$
\ecor


In \cite{oA}, the authors got the following homogeneous expansion of $\alpha$-harmonic functions
(see \cite[Theorem 1.2]{oA}):

A function $f$ in $\mathbb{D}$ is $\alpha$-harmonic if and only if it has the following convergent power series expansion:
\beq\label{eq1.6}\;\;\;
f(z)=\sum^{\infty}_{k=0}c_{k}z^{k}+\sum^{\infty}_{k=1}c_{-k}P_{\alpha,k}(|z|^2)\overline{z}^{k},
\eeq
where $P_{\alpha,k}(x)=\int^{1}_{0}t^{k-1}(1-tx)^{\alpha}dt$ $(-1<x<1)$ and $\{c_{k}\}^{\infty}_{k=-\infty}$ denotes a sequence of complex numbers with
$
\lim_{|k|\to\infty}\sup|c_{k}|^{\frac{1}{|k|}}\leq1.
$

The second aim of this paper is to prove the following estimates on coefficients $c_{k}$ and $c_{-k}$.

\begin{thm}\label{thm1.3}
Suppose that $f$ is an $\alpha$-harmonic function in $\ID$ with $\alpha>-1$ and that
$\sup_{z\in \mathbb{D}}|f(z)|\leq M$, where $M$ is a constant. If $f$ has the series expansion \eqref{eq1.6}, then for $k\in\{0,1,2,...\},$
\be\label{sun-1}|c_{k}|\leq M,\ee
and for $k\in\{1,2,...\},$
\be\label{sun-2}|c_{k}|+|c_{-k}|B(k, \alpha+1)\leq\frac{4M}{\pi},\ee
 where $B(p, q)$ denotes the Beta function.
\end{thm}

By \cite[Definition 2.1]{oA}, we find that
$$\mathcal{P}_{\alpha}(ze^{-i\theta})=\sum^{\infty}_{k=0}e^{-ik\theta}z^{k}
+\sum^{\infty}_{k=1}\frac{\Gamma(k+\alpha+1)}{\Gamma(k)\Gamma(\alpha+1)}P_{\alpha,k}(|z|^2)e^{ik\theta}\overline{z}^{k}.$$
If $|f^{\ast}(z)|\leq M$, then by \eqref{eq1.4}, we get
$$|c_{-k}|=\left|\frac{\Gamma(k+\alpha+1)}{\Gamma(k)\Gamma(\alpha+1)}\frac{1}{2\pi}\int^{2\pi}_{0}e^{ik\theta}f^{\ast}(e^{i\theta})d\theta\right|
\leq M\frac{\Gamma(k+\alpha+1)}{\Gamma(k)\Gamma(\alpha+1)}\rightarrow \infty$$
as $k\rightarrow\infty$.

Moreover, from the proof of \cite[Theorem 1.2]{oA}, we see that
$$(1-{|{z}|}^{2})^{-\alpha}\frac{\partial}{\partial \overline{z}}f(z)=\overline{h(z)},$$
where $h(z)=\sum^{\infty}_{k=0}a_{k}z^k$, $z\in\ID$ and $c_{-k}=\overline{a_{k-1}}$ for $k\geq1$. Note that if $h(z)$ is a normalized (in the sense that $h(0)=h'(0)-1=0$) univalent analytic function in $\ID$, then by Louis de Branges's theorem
it is well-known that $|a_{k}|\leq k$ for all $k\geq2$ so that
\be\label{re-1}c_{-1}=0,\;\; c_{-2}=1\;\; \mbox{and}\;\; |c_{-k}|=|a_{k-1}|\leq k-1\;\;\mbox{for all}\;\;k\geq 3.\ee

The classical Landau theorem says that there is a $\rho=\frac{1}{M+\sqrt{M^{2}-1}}$ such that every function $f$, analytic in $\mathbb{D}$ with $f(0)=f^{'}(0)-1=0$ and $|f(z)|< M$, is univalent in the disk $\mathbb{D}_{\rho}$. Moreover, the range $f(\mathbb{D}_{\rho})$ contains a disk of radius $M\rho^{2}$, where $M\geq1$ is a constant (see \cite{E}). Recently, many authors considered Landau type theorem for $\alpha$-harmonic functions $f$ when $\alpha=0$ (see \cite{CPW, ChW, CSP, cSW, cpx, CSA} etc).\medskip

As an application of Theorems \ref{thm1.1} and \ref{thm1.3}, we get the following Landau type theorem for $\alpha$-harmonic functions.

\begin{thm}\label{thm1.4}
Suppose that $f$ is an $\alpha$-harmonic function in $\ID$ with $\alpha\geq 0$, that $f^{*}\in C(\mathbb{T})$,
that $\sup_{z\in\mathbb{\overline{D}}}|f(z)|\leq M$, where $M$ is a constant, and that $f(0)=|J_{f}(0)|-\beta=0$. If $f$ satisfies \eqref{re-1},
 then we have the following:
\begin{enumerate}
  \item \label{eq1.7}
  $f$ is univalent in $\mathbb{D}_{\rho_{0}}$, where $\rho_{0}$ satisfies the following equation
\be\label{fil-1}
\frac{\beta\mathrm{c}_{\alpha}}{M(\alpha+2)}-(M+5)\,\frac{\rho_{0}(2-\rho_{0})}{(1-\rho_{0})^{2}}=0;
\ee
\item \label{eq1.10}
  $f(\mathbb{D}_{\rho_{0}})$ contains a univalent disk $\mathbb{D}_{R_{0}}$ with
$$
R_{0}\geq(M+5)\,\Big(\frac{\rho_{0}}{1-\rho_{0}}\Big)^2.
$$
\end{enumerate}
\end{thm}

The arrangement of the rest of this paper is as follows. In Section \ref{sec-2}, we shall prove Theorem \ref{thm1.1} and Corollary \ref{cor1.1}. Section \ref{sec-3} will be devoted to the proof of Theorem \ref{thm1.3}. In Section \ref{sec-4}, Theorem \ref{thm1.4} will be demonstrated.

\section{Schwarz-Pick type inequality}\label{sec-2}

The aim of this section is to prove Theorem \ref{thm1.1} and Corollary \ref{cor1.1}.
The proofs need a result from \cite{OD}. Before the statement of this result, we do some preparation.

In \cite{OD}, the author considered the following integral means:
\beq\label{eq1.5}\;\;\;
\mathcal{M}_{\alpha}(r)= \frac{1}{2\pi} \int_{0}^{2\pi}K_{\alpha}(re^{i\theta})d\theta,
\eeq
where $r\in[0, 1)$ and
$$
K_{\alpha}(z)=\mathrm{c}_{\alpha}|\mathcal{P}_{\alpha}(z)|=\mathrm{c}_{\alpha}\frac{(1-|z|^{2})^{\alpha+1}}{|1-z|^{\alpha+2}}
$$ in $\ID$.

Let us recall the following result from \cite{OD}.

\begin{Thm}\label{Thm A} \cite[Theorem 3.1]{OD}
Let $\alpha>-1$. The integral means function $\mathcal{M}_{\alpha}(r)$ given by \eqref{eq1.5}
satisfies the following assertions.
\ben
\item
$\lim_{r\to1^{-}}\mathcal{M}_{\alpha}(r)=1;$
\item
$\mathcal{M}_{\alpha}^{(n)}(r)\geq0$ for $r\in[0,1)$ and $n\geq0$.
\een
\end{Thm}

The following result also plays a key role in the proof of Theorem \ref{thm1.1}.

\begin{lem}\label{lem1}
If $\alpha>-1$ and $f^{*}\in C(\mathbb{T})$, then
$$\frac{\partial}{\partial z}\int^{2\pi}_{0}\mathcal{P}_{\alpha}(ze^{-i\theta})f^{\ast}(e^{i\theta})d\theta
=\int^{2\pi}_{0}\frac{\partial}{\partial z}\mathcal{P}_{\alpha}(ze^{-i\theta})f^{\ast}(e^{i\theta})d\theta$$
and
$$\frac{\partial}{\partial \overline{z}}\int^{2\pi}_{0}\mathcal{P}_{\alpha}(ze^{-i\theta})f^{\ast}(e^{i\theta})d\theta
=\int^{2\pi}_{0}\frac{\partial}{\partial \overline{z}}\mathcal{P}_{\alpha}(ze^{-i\theta})f^{\ast}(e^{i\theta})d\theta.$$
\end{lem}

\bpf
By elementary calculations we see that the following equalities hold:
\beq\label{eq11}\frac{\partial}{\partial{z}}\mathcal{P}_{\alpha}(ze^{-i\theta})
=\frac{(1-|z|^{2})^{\alpha}\left[e^{-i\theta}(1-|z|^{2})-(\alpha+1)\overline{z}(1-ze^{-i\theta})\right]}
{(1-ze^{-i\theta})^{2}(1-\overline{z}e^{i\theta})^{\alpha+1}}\eeq
and
\beq\label{eq12}
\frac{\partial}{\partial{\overline{z}}}\mathcal{P}_{\alpha}(ze^{-i\theta})=\frac{(\alpha+1)(1-|z|^{2})^{\alpha}e^{i\theta}}{(1-\overline{z}e^{i\theta})^{\alpha+2}}.
\eeq
Then we know that functions
$$\frac{\partial}{\partial z}\mathcal{P}_{\alpha}(ze^{-i\theta})f^{\ast}(e^{i\theta})\;\;\mbox{and}\;\;
\frac{\partial}{\partial \overline{z}}\mathcal{P}_{\alpha}(ze^{-i\theta})f^{\ast}(e^{i\theta})$$
are continuous on $\overline{\ID}_r\times[0, 2\pi]$, where $r\in[0, 1)$.

Let $z=\rho e^{i\varphi}\in \overline{\ID}_r$. It follows from
$$\frac{\partial }{\partial \rho} \mathcal{P}_{\alpha}(ze^{-i\theta})=\frac{\partial }{\partial z} \mathcal{P}_{\alpha}(ze^{-i\theta})e^{i\varphi}
+\frac{\partial }{\partial \overline{z}}\mathcal{P}_{\alpha}(ze^{-i\theta})e^{-i\varphi}$$
and
$$\frac{\partial }{\partial \varphi}\mathcal{P}_{\alpha}(ze^{-i\theta})=\frac{\partial }{\partial z} \mathcal{P}_{\alpha}(ze^{-i\theta})iz
-\frac{\partial }{\partial \overline{z}}\mathcal{P}_{\alpha}(ze^{-i\theta})i\overline{z}$$
that both
$$\frac{\partial }{\partial \rho} \mathcal{P}_{\alpha}(ze^{-i\theta})f(e^{i\theta})\;\;\mbox{ and} \;\;\frac{\partial }{\partial \varphi}\mathcal{P}_{\alpha}(ze^{-i\theta})f(e^{i\theta})$$
are continuous in $\overline{\ID}_r\times[0, 2\pi]$.
Hence
\beqq
\int^{\rho}_{0}\int^{2\pi}_{0}\frac{\partial }{\partial \rho} \mathcal{P}_{\alpha}(ze^{-i\theta})f^{\ast}(e^{i\theta})d\theta d\rho
&=&\int^{2\pi}_{0}\int^{\rho}_{0}\frac{\partial }{\partial \rho} \mathcal{P}_{\alpha}(ze^{-i\theta})f^{\ast}(e^{i\theta})d\rho d\theta\\
&=&\int^{2\pi}_{0}\big(\mathcal{P}_{\alpha}(ze^{-i\theta})-\mathcal{P}_{\alpha}(0)\big)f^{\ast}(e^{i\theta})d\theta
\eeqq
and
\beqq
\int^{\varphi}_{0}\int^{2\pi}_{0}\frac{\partial }{\partial \varphi}\mathcal{P}_{\alpha}(ze^{-i\theta})f^{\ast}(e^{i\theta})d\theta d\varphi
&=&\int^{2\pi}_{0}\int^{\varphi}_{0}\frac{\partial }{\partial \varphi}\mathcal{P}_{\alpha}(ze^{-i\theta})f^{\ast}(e^{i\theta})d\varphi d\theta\\
&=&\int^{2\pi}_{0}\big(\mathcal{P}_{\alpha}(ze^{-i\theta})-\mathcal{P}_{\alpha}(\rho e^{-i\theta})\big)f^{\ast}(e^{i\theta})d\theta.
\eeqq

By differentiating with respect to $\rho$ and $\varphi$, respectively, we get
\beq\label{eq13}
\int^{2\pi}_{0}\frac{\partial }{\partial \rho} \mathcal{P}_{\alpha}(ze^{-i\theta})f^{\ast}(e^{i\theta})d\theta=\frac{\partial}{\partial \rho}\int^{2\pi}_{0}\mathcal{P}_{\alpha}(ze^{-i\theta})f^{\ast}(e^{i\theta})d\theta
\eeq
and
\beq\label{eq14}
\int^{2\pi}_{0}\frac{\partial }{\partial \varphi}\mathcal{P}_{\alpha}(ze^{-i\theta})f^{\ast}(e^{i\theta})d\theta=\frac{\partial}{\partial \varphi}\int^{2\pi}_{0}\mathcal{P}_{\alpha}(ze^{-i\theta})f^{\ast}(e^{i\theta})d\theta.
\eeq
Since
$$\frac{\partial }{\partial z} \mathcal{P}_{\alpha}(ze^{-i\theta})=\frac{e^{-i\varphi}}{2}\left(\frac{\partial }{\partial \rho} \mathcal{P}_{\alpha}(ze^{-i\theta})-\frac{i}{\rho}\frac{\partial }{\partial \varphi}\mathcal{P}_{\alpha}(ze^{-i\theta})\right)$$
and
$$\frac{\partial }{\partial \overline{z}}\mathcal{P}_{\alpha}(ze^{-i\theta})=\frac{e^{i\varphi}}{2}\left(\frac{\partial }{\partial \rho} \mathcal{P}_{\alpha}(ze^{-i\theta})+\frac{i}{\rho}\frac{\partial }{\partial \varphi}\mathcal{P}_{\alpha}(ze^{-i\theta})\right),$$
it follows from \eqref{eq13} and \eqref{eq14} that the proof of the lemma is complete.
\epf

Now, we are ready to present the proofs of Theorem \ref{thm1.1} and Corollary \ref{cor1.1}.\medskip

\noindent {\bf Proof of Theorem \ref{thm1.1}}
From \eqref{eq11} and \eqref{eq12}, we can easily get
$$
\left|\frac{\partial}{\partial{z}}\mathcal{P}_{\alpha}(ze^{-i\theta})\right|\leq \frac{1}{\mathrm{c}_{\alpha}}\cdot\frac{(\alpha+2)|z|+1}{1-|z|^{2}}K_{\alpha}(ze^{-i\theta})
$$
and
$$\left|\frac{\partial}{\partial{\overline{z}}}\mathcal{P}_{\alpha}(ze^{-i\theta})\right|= \frac{\alpha+1}{\mathrm{c}_{\alpha}}\cdot\frac{1}{1-|z|^{2}}K_{\alpha}(ze^{-i\theta}).
$$
In the first inequality above, the fact ``$1-|z|\leq |1-ze^{-i\theta}|$" is applied. By \eqref{qh-1}, \eqref{eq1.4} and Lemma \ref{lem1} yield
$$|D{f}(z)|=\left|\frac{1}{2\pi}\int_{0}^{2\pi}\frac{\partial}{\partial{z}}\mathcal{P}_{\alpha}(ze^{-i\theta})f^{\ast}(e^{i\theta})d\theta\right|
+\left|\frac{1}{2\pi}\int_{0}^{2\pi}\frac{\partial}{\partial{\overline{z}}}\mathcal{P}_{\alpha}(ze^{-i\theta})f^{\ast}(e^{i\theta})d\theta\right|,
$$
we see from \eqref{eq1.5} and Theorem \Ref{Thm A} that
\begin{eqnarray*}
|D{f}(z)|
\leq\frac{M(\alpha+2)}{\mathrm{c}_{\alpha}}\cdot\frac{1}{1-|z|}\mathcal{M}_{\alpha}(|z|)
\leq \frac{M(\alpha+2)}{\mathrm{c}_{\alpha}}\cdot\frac{1}{1-|z|},
\end{eqnarray*}
and so the proof of Theorem \ref{thm1.1} is complete. \qed
\medskip

\noindent {\bf Proof of Corollary \ref{cor1.1}} For any $z_1$ and $z_2\in \ID$, let $\gamma$ be the hyperbolic geodesic connecting $z_1$ and $z_2$. It follows from Theorem \ref{thm1.1} that
$$
|f(z_1)-f(z_2)|\leq\int_{\gamma}|D{f}(z)|\cdot|dz|
\leq \frac{\alpha+2}{\mathrm{c}_{\alpha}}\int_{\gamma}\frac{2}{1-|z|^2}|dz|
= \frac{\alpha+2}{\mathrm{c}_{\alpha}}d_{h_\mathbb{D}}(z_1, z_2),
$$
as required. \qed

\section{Estimates on coefficients}\label{sec-3}

The aim of this paper is to prove Theorem \ref{thm1.3}. We start with a lemma.

\blem\label{sun-4}
Under the assumptions of Theorem \ref{thm1.3}, if $f$ has the series expansion \eqref{eq1.6}, then \ben
\item\label{sat-1}
$|c_k|\leq  M$ for $k\geq0$;
\item\label{sat-2}
$(|c_{k}|+|c_{-k}|P_{\alpha,k}(r^{2}))r^{k}\leq \frac{4}{\pi}M$ for $k>0$ and $r\in (0, 1)$.
\een
\elem
\bpf If $k\not=0$, let $z=re^{i\theta}\in{\mathbb{D}}$. Then by \eqref{eq1.6}, we have
$$
c_{k}r^{k}=\frac{1}{2\pi}\int_{0}^{2\pi}f(re^{i\theta})e^{-ik\theta}d{\theta}\;\;\mbox{and}\;\;
c_{-k}P_{\alpha,k}(r^{2})r^{k}=\frac{1}{2\pi}\int_{0}^{2\pi}f(re^{i\theta})e^{ik\theta}d{\theta}.
$$

Letting $c_k=|c_k|e^{i\mu_k}$ and  $c_{-k}=|c_{-k}|e^{i\nu_k}$ leads to
\beq\label{eq2.6}\nonumber
(|c_{k}|+|c_{-k}|P_{\alpha,k}(r^{2}))r^{k}
&=&\frac{1}{2\pi}\int^{2\pi}_{0}f(re^{i\theta})\left(e^{-i(k\theta+\mu_k)}+e^{i(k\theta-\nu_k)}\right)d{\theta}
\\\nonumber
&\leq&\frac{1}{2\pi}\int^{2\pi}_{0}|f(re^{i\theta})|\cdot\left|e^{-i(k\theta+\mu_k)}+e^{i(k\theta-\nu_k)}\right|d{\theta}
\\\nonumber
&\leq&\frac{M}{\pi}\int^{2\pi}_{0}\left|\cos\left(k\theta+\frac{\mu_k-\nu_k}{2}\right)\right|d{\theta},
\eeq
and so
\cite[Lemma 1]{cR} gives
$$
(|c_{k}|+|c_{-k}|P_{\alpha,k}(r^{2}))r^{k}
\leq \frac{4M}{\pi}.
$$
Thus the assertion \eqref{sat-2} in the lemma is true.

To prove the assertions \eqref{sat-1}, we first recall from \cite[Definition 2.1]{oA} that
$$\mathcal{P}_{\alpha}(ze^{-i\theta})=\sum^{\infty}_{k=0}e^{-ik\theta}z^{k}
+\sum^{\infty}_{k=1}\frac{\Gamma(k+\alpha+1)}{\Gamma(k)\Gamma(\alpha+1)}P_{\alpha,k}(|z|^2)e^{ik\theta}\overline{z}^{k}.$$
Then by \eqref{eq1.4}, we get
\beqq
f(z)&=&\sum^{\infty}_{k=0}z^k\frac{1}{2\pi}\int^{2\pi}_{0}e^{-ik\theta}f^{\ast}(e^{i\theta})d{\theta}\\
&&+\sum^{\infty}_{k=1}\frac{\Gamma(k+\alpha+1)}{\Gamma(k)\Gamma(\alpha+1)}P_{\alpha,k}(|z|^2)\overline{z}^{k}
\frac{1}{2\pi}\int^{2\pi}_{0}e^{ik\theta}f^{\ast}(e^{i\theta})d{\theta},
\eeqq
which implies
$$
|c_k|=\left|\frac{1}{2\pi}\int^{2\pi}_{0}e^{-ik\theta}f^{\ast}(e^{i\theta})d{\theta}\right|\leq M,
$$
as required.
\epf

\noindent {\bf Proof of Theorem \ref{thm1.3}}\; To prove this theorem, by Lemma \ref{sun-4}, we only need to check \eqref{sun-2} in the theorem.
By letting $r\to1^{-}$ in Lemma \ref{sun-4}\eqref{sat-2}, we see that the inequalities
\eqref{sun-2} easily follows.
\qed

\section{Landau type theorem }\label{sec-4}
This section consists of two subsections. In the first subsection, we shall prove an auxiliary result. In the second subsection, Theorem \ref{thm1.4} will be checked.

\subsection{A lemma}
\begin{lem}\label{lemm3.1}
For constants $\alpha>-2$, $\beta>0$ and $M>0$, let
$$
\varphi(x)=\frac{\beta\mathrm{c}_{\alpha}}{M(\alpha+2)}+(M+5)\,\frac{x(x-2)}{(1-x)^{2}}
$$
in $[0,1)$. Then \ben
\item
 $\varphi$ is continuous in $[0,1)$  and strictly decreasing in $(0,1)$;
 \item
 there is a unique $x_{0}\in (0,1)$ such that $\varphi(x_{0})=0$.
 \een
\end{lem}
\bpf
For $x\in[0,1)$, obviously,
$$
\varphi'(x)=-\frac{2(M+5)}{(1-x)^{3}}<0.
$$
Hence $\varphi(x)$ is continuous and strictly decreasing in $[0,1)$.
It follows from
$$
\varphi(0)=\frac{\beta\mathrm{c}_{\alpha}}{M(\alpha+2)}>0 \;\;\text{and}\;\; \lim_{x\to1^{-}}\varphi(x)=-\infty<0
$$
that there is a unique $x_{0}\in (0,1)$ such that $\varphi(x_{0})=0$. The proof of this lemma is complete.
\epf

\subsection{Proof of Theorem \ref{thm1.4}.}
To prove this theorem, we need estimates on two quantities $|f_{z}(z)-f_{z}(0)|+|f_{\overline{z}}(z)-f_{\overline{z}}(0)|$ and $l(D{f}(0))$.
First, we estimate $|f_{z}(z)-f_{z}(0)|+|f_{\overline{z}}(z)-f_{\overline{z}}(0)|$. Obviously,
by \eqref{eq1.6}, we see that
$$
f_{z}(z)-f_{z}(0)=\sum_{k=2}^{\infty}kc_{k}z^{k-1}+\sum_{k=2}^{\infty}c_{-k}\frac{d}{dw}P_{\alpha, k}(w)\overline{z}^{k+1}
$$
and
$$
f_{\overline{z}}(z)-f_{\overline{z}}(0)=\sum_{k=2}^{\infty}kc_{-k}P_{\alpha, k}(w)\overline{z}^{k-1}+\sum_{k=2}^{\infty}c_{-k}\frac{d}{dw}P_{\alpha, k}(w)z\overline{z}^{k},
$$
where $w=|z|^{2}$.

Since
$$\frac{d}{dw}P_{\alpha,k}(w)=-\int^{1}_{0}t^{k}\alpha(1-tw)^{\alpha-1}dt\leq 0,$$
we get that
\beq\label{eq3.1.0}\;\;\;
P_{\alpha,k}(w)\leq P_{\alpha,k}(0)=\frac{1}{k}.
\eeq

Moreover, since
$$
P_{\alpha, k}(w)=\frac{1}{w^{k}}\int^{w}_{0}x^{k-1}(1-x)^{\alpha}dx,
$$ we easily get
\beq\label{eq3.1}\;\;\;
\frac{d}{dw}P_{\alpha,k}(w)=-\frac{k}{w}P_{\alpha,k}(w)+\frac{(1-w)^\alpha}{w}.
\eeq

Then \eqref{sun-2}, \eqref{re-1}, \eqref{eq3.1.0} and \eqref{eq3.1} guarantee that
\beq\label{eq3.2}\;\;\;
|f_{z}(z)-f_{z}(0)|+|f_{\overline{z}}(z)-f_{\overline{z}}(0)|
&\leq&\sum_{k=2}^{\infty}k\big(|c_{k}|+|c_{-k}|P_{\alpha,k}(w)\big)|z|^{k-1}\\\nonumber &&
+2\sum_{k=2}^{\infty}|c_{-k}|\big(kP_{\alpha,k}(w)+1\big)|z|^{k-1}
\\\nonumber
&\leq&(M+5)\sum_{k=2}^{\infty}k|z|^{k-1}
\\\nonumber
&=&(M+5)\,\frac{|z|(2-|z|)}{(1-|z|)^{2}},
\eeq which is what we want.

Next, we estimate $l(D{f}(0))$. Applying Theorem \ref{thm1.1} leads to
$$
\beta=|J_{f}(0)|=|D{f}(0)|\;l(D{f}(0))\leq \frac{M(\alpha+2)}{\mathrm{c}_{\alpha}}\;l(D{f}(0)),
$$
which gives
\beq\label{eq3.3}\;\;\;
l(D{f}(0))\geq\frac{\beta\mathrm{c}_{\alpha}}{M(\alpha+2)}.
\eeq

Now, we are ready to finish the proof of the theorem. First, we demonstrate the univalence of $f$ in $\mathbb{D}_{\rho_{0}}$, where $\rho_{0}$ is determined by the equation \eqref{fil-1}. For this, let $z_{1}$, $z_{2}$ be two points in $\mathbb{D}_{\rho_{0}}$ with $z_{1}\neq z_{2}$, and denote the segment from $z_{1}$ to $z_{2}$ with the endpoints $z_{1}$ and $z_{2}$ by $[z_{1},z_{2}]$. Since
\begin{eqnarray*}
|f(z_{2})-f(z_{1})|&=&\left|\int_{[z_{1},z_{2}]}f_{z}(z)dz+f_{\overline{z}}(z)d\overline{z}\right|
\\\nonumber
&\geq&\left|\int_{[z_{1},z_{2}]}f_{z}(0)dz+f_{\overline{z}}(0)d\overline{z}\right|\\\nonumber
&&
-\left|\int_{[z_{1},z_{2}]}[f_{z}(z)-f_{z}(0)]dz+[f_{\overline{z}}(z)-f_{\overline{z}}(0)]d\overline{z}\right|
,\end{eqnarray*}
we see from \eqref{eq3.2}, \eqref{eq3.3} and Lemma \ref{lemm3.1} that

\begin{eqnarray*}
|f(z_{2})-f(z_{1})|
&\geq& l(D{f}(0))\cdot|z_{2}-z_{1}|-(M+5)\int_{0}^{|z_{2}-z_{1}|}\frac{|z|(2-|z|)}{(1-|z|)^{2}}|dz|
\\\nonumber&>&\left[\frac{\beta\mathrm{c}_{\alpha}}{M(\alpha+2)}-(M+5)\,\frac{\rho_{0}(2-\rho_{0})}{(1-\rho_{0})^{2}}\right]|z_{2}-z_{1}|
\\\nonumber
&=&0.
\end{eqnarray*}
Thus, for arbitrary $z_{1}$ and $z_{2}\in\mathbb{D}_{\rho_{0}}$ with $z_{1}\neq z_{2}$, we have $$f(z_{1})\neq f(z_{2}),$$ which implies the univalence of $f$ in $\mathbb{D}_{\rho_{0}}$.

Next, we prove Theorem \ref{thm1.4}\eqref{eq1.10}. For any $\zeta=\rho_{0}e^{i\theta}\in \partial\mathbb{D}_{\rho_{0}}$, we obtain that

\begin{eqnarray*}
|f(\zeta)-f(0)|&=&\left|\int_{[0,\zeta]}f_{z}(z)dz+f_{\overline{z}}(z)d\overline{z}\right|
\\\nonumber
&\geq&\left|\int_{[0,\zeta]}f_{z}(0)dz+f_{\overline{z}}(0)d\overline{z}\right|\\ \nonumber
&&-\left|\int_{[0,\zeta]}[f_{z}(z)-f_{z}(0)]dz+[f_{\overline{z}}(z)-f_{\overline{z}}(0)]d\overline{z}\right|
\\\nonumber
&\geq&l(D{f}(0))\rho_{0}-(M+5)\int_{0}^{\rho_0}\frac{|z|(2-|z|)}{(1-|z|)^{2}}|dz|\hspace{1cm}(\text{by \eqref{eq3.2}})
\\\nonumber
&=& \frac{\beta\mathrm{c}_{\alpha}\rho_{0}}{M(\alpha+2)}-(M+5)\,\frac{\rho_{0}^2}{1-\rho_{0}}.
\\\nonumber
&=& (M+5)\,\Big(\frac{\rho_{0}}{1-\rho_{0}}\Big)^2. \hspace{1cm}(\text{by \eqref{fil-1}})
\end{eqnarray*}
Hence $f(\mathbb{D}_{\rho_{0}})$ contains a univalent disk $\mathbb{D}_{R_{0}}$, where
$$
R_{0}\geq(M+5)\,\Big(\frac{\rho_{0}}{1-\rho_{0}}\Big)^2.
$$
The proof of this theorem is complete.\qed
\bigskip\bigskip

\begin{center}
{\bf Acknowledgements}
\end{center}

The research was partly supported by the NSFs of China (No. 11571216 and No. 11671127), the Hunan Provincial Innovation Foundation For Postgraduate (No. CX2016B159), the NSF of Guangdong Province (No. 2014A030313471) and Project of ISTCIPU in Guangdong Province (No. 2014KGJHZ007).

The authors thank the referee very much for his/her careful reading of this paper and many useful suggestions.


\bigskip


\begin{thebibliography}{99}

\bibitem{Bea} A. F.
Beardon, The geometry of discrete groups, Graduate Texts in Mathematics, 91. Springer-Verlag, New York, 1983.


\bibitem{HH} H. Chen,
The Schwarz-Pick lemma for planar harmonic mappings,
\textit{Sci. China Math.}, {\bf 54} (2011), 1101--1118.

\bibitem{CPW} H. Chen, P. Gauthier and W. Hengartner,
Bloch constants for planar harmonic mappings,
\textit{Proc. Amer. Math. Soc.}, {\bf 128} (2000), 3231--3240.

\bibitem{cR} Sh. Chen, A. Rasila,
Schwarz-Pick type estimates of pluriharmonic mappings in the unit polydisk,
\textit{Illinois J. Math.}, {\bf 58} (2014), 1015--1024.

\bibitem{ChW}Sh. Chen, S. Ponnusamy and X. Wang,
Integral means and coefficient estimates on planar harmonic mappings,
\textit{Ann. Acad. Sci. Fenn. Math.}, {\bf 37} (2012), 69--79.

\bibitem{CSP} Sh. Chen, S. Ponnusamy and X. Wang,
On planar harmonic Lipschitz and planar harmonic Hardy classes,
\textit{Ann. Acad. Sci. Fenn. Math.}, {\bf 36} (2011), 567--576.

\bibitem{cSW} Sh. Chen, S. Ponnusamy and X. Wang,
Harmonic mappings in Bergman spaces,
\textit{Monatsh. Math.}, {\bf 170} (2013), 325--342.

\bibitem{cpx} Sh. Chen, S. Ponnusamy and X. Wang,
Bloch constant and Landau's theorem for planar p-harmonic mappings,
\textit{J. Math. Anal. Appl.}, {\bf 373} (2011), 102--110.

\bibitem{CSA}Sh. Chen, S. Ponnusamy and A. Rasila,
Coefficient estimates, Landau¡¯s theorem and Lipschitz-type spaces on planar harmonic mappings,
\textit{J. Aust. Math. Soc.}, {\bf 96} (2014), 198--215.

\bibitem{CM} Sh. Chen, M. Vuorinen,
Some properties of a class of elliptic partial differential operators.
\textit{J. Math. Anal. Appl.}, {\bf431} (2015), 1124--1137.

\bibitem{F} F. Colonna,
The Bloch constant of bounded harmonic mappings,
\textit{Indiana Univ. Math. J.}, {\bf 38} (1989), 829--840.

\bibitem{DuR} P. Durin,
Harmonic mappings in the plane,
Cambridge Tracts in Mathematics, 156. Cambridge University Press, Cambridge, 2004.

\bibitem{HKZ} H. Hedenmalm, B. Korenblum and K. Zhu,
Theory of Bergman spaces,
\textit{Springer-Verlag,} 2000.

\bibitem{HeO} H. Hedenmalm and A. Olofsson, Hele-Shaw flow on weakly hyperbolic surfaces,
\textit{Indiana Univ. Math. J.,} {\bf 54} (2005), 1161--1180.

\bibitem{HeP} H. Hedenmalm and Y. Perdomo, Mean value surfaces with prescribed curvature form,
\textit{J. Math. Pures Appl.,} {\bf 83} (2004), 1075--1107.

\bibitem{HeS} H. Hedenmalm and S. Shimorin, Hele-Shaw flow on hyperbolic surfaces,
\textit{J. Math. Pures Appl.,} {\bf 81} (2002), 187--222.



\bibitem{KaM} D. Kalaj and M. Vuorinen,
On harmonic functions and the Schwarz lemma,
\textit{Proc. Amer. Math. Soc.}, {\bf140} (2012), 161--165.

\bibitem{E} E. Landau,
\"Uber die Bloch¡¯sche konstante und zwei verwandte weltkonstanten,
\textit{ Math. Zeit.}, {\bf 30} (1929), 608--634.

\bibitem{OD} A. Olofsson,
Differential operators for a scale of Poisson type kernels in the unit disc,
\textit{J. Anal. Math.}, {\bf 123} (2014), 227--249.

\bibitem{oA} A. Olofsson and J. Wittsten,
Poisson integrals for standard weighted Laplacians in the unit disc,
\textit{J. Math. Soc. Japan}, {\bf 65} (2013), 447--486.

\bibitem{M} M. Pavlovi\`{c},
Harmonic Schwarz lemmas: Chen, Kalaj-Vuorinen, Pavlovi\`{c} and Heinz, preprint.

\bibitem{SSh} S. Shimorin, On Beurling-type theorems in weighted $\ell^2$ and Bergman spaces, \textit{Proc. Amer.
Math. Soc.}, {\bf 131} (2003), 1777--1787.


\end{thebibliography}
\end{document}